\documentstyle[amssymb,amsfonts]{amsart}

\def\be#1{\begin{equation} \label{#1}}
\def\bs{\begin{split}}
\def\es{\end{split}}
\def\ba{\begin{align}}
\def\bas{\begin{align*}}
\def\R{{\mbox{\bf R}}}

\def\eps{\varepsilon}

\def\emph#1{{\it #1}}
\def\textbf#1{{\bf #1}}
\newenvironment{proof}{\noindent {\bf Proof} }{\endprf\par}
\def \endprf{\hfill  {\vrule height6pt width6pt depth0pt}\medskip}

\theoremstyle{plain}
  \newtheorem{theorem}[subsection]{Theorem}
  
  \newtheorem{lemma}[subsection]{Lemma}
  \newtheorem{corollary}[subsection]{Corollary}
  \newtheorem{conjecture}[subsection]{Conjecture}

\theoremstyle{remark}

\theoremstyle{definition}

\include{psfig}

\begin{document}

\title[Ill-posedness for one-dimensional wave maps]{
Ill-posedness for one-dimensional wave maps at the critical regularity}

\author{Terence Tao}
\address{Department of Mathematics, UCLA, Los Angeles, CA 90024}
\email{tao@@math.ucla.edu}

\begin{abstract}
We show that the wave map equation in $\R^{1+1}$ is in general ill-posed
in the critical space $\dot H^{1/2}$, and the Besov space
$\dot B^{1/2,1}_2$.  The problem is attributed to the bad behaviour
of the one-dimensional bilinear expression $D^{-1}(f Dg)$ in
these spaces.
\end{abstract}

\maketitle

\section{Introduction}\label{intro}

Write $(\R^{n+1}, g)$ for $n+1$ dimensional  Minkowski space
with flat metric $g = \text{diag}(-1,1,\ldots 1)$. In what follows
$(M,h)$ will denote a Riemannian manifold with metric $h$.

We consider the Cauchy problem
\be{wavemap}
\bs
\Box \phi^k+
{\Gamma}_{\alpha \beta}^k(\phi) \partial_\mu\phi^\alpha\partial^{\mu}\phi^{\beta
}
& = 0.
\\
\phi(0) &= f\\
\phi_t(0) &= 0
\es
\end{equation}
 in local co-ordinates
where $\phi: \R^{n+1} \to M$ is a map, $\Gamma^{k}_{\alpha \beta}$
are the Christoffel symbols
corresponding to the Riemannian metric $h$, and the initial data $f$
takes values in $M$.  To avoid technical issues involving negative order
Sobolev spaces, we will always take the initial velocity $\phi_t$ to be zero.

Solutions to \eqref{wavemap} are called wave maps with initial data $f$.
For a discussion of the significance of these maps and some open problems
we refer to \cite{struwe.barrett},\cite{kman.barrett},\cite{grillakis.zurich},
and \cite{shatah.zurich}.  For this
paper we restrict ourselves to the case when the target manifold is a sphere
$M = S^{m-1}$, which we imbed into Euclidean space in the usual manner $\R^m$.
The problem then simplifies to
\be{wave-map}
\bs
\Box \phi+
\phi (\partial_\mu\phi \cdot \partial^{\mu}\phi)
& = 0, \\
\phi(0) &= f\\
\phi_t(0) &= 0\\
f(x) &\in S^{m-1}.
\es
\end{equation}
See e.g. \cite{struwe.barrett}.
One can show that solutions to \eqref{wave-map} stay on the sphere, at least
if $\phi$ is smooth.

We consider the problem of whether this Cauchy problem is locally well-posed for 
small data $f$ in various spaces of Sobolev type.  It is known that the problem
is locally well-posed in $H^s$ for $s > n/2$ (see \cite{kman.mach.smoothing},
\cite{kman.selberg}).
In the case $n=1$ then one also has global well-posedness if $s > 3/4$ (\cite{keel:wavemap}; see also
earlier results in \cite{pohl}, \cite{gu}, \cite{shatah}, \cite{lady}).
When $s < n/2$ the problem is super-critical and the problem is almost certainly
ill-posed for the sphere.  (For instances of blowup in the supercritical case,
see \cite{shatah}, \cite{shatah.shadi.blow}).  The critical case $s = n/2$ 
appears to be very difficult, and only partial results are known.  In the critical case
the problem is invariant under scaling, and one should replace the inhomogeneous
Sobolev spaces $H^s$ with the homogeneous spaces $\dot H^s$.  We will adopt the convention
that the constant functions have zero $\dot H^s$ norm, since it would otherwise be impossible
to have functions taking values in $S^{m-1}$ which have finite $\dot H^s$ norm.  In practice
we will subtract a constant (such as the unit vector $e_1$) from functions on $S^{n-1}$
in order to measure such norms.

\begin{conjecture}\label{hs-c}  The Cauchy problem \eqref{wave-map} is well-posed for small data
in $\dot H^{\frac{n}{2}}$.
\end{conjecture}

Note that by scale invariance considerations, the problem of local well-posedness is equivalent
to that of global well-posedness in the small-data critical case.  We note that even if
the above Cauchy problem is well-posed, the solution operator cannot be twice
differentiable in $\dot H^{\frac{n}{2}}$; 
see \cite{keel:wavemap}.  This indicates that the problem cannot be
attacked by purely iterative methods, and requires some additional subtlety.

Recently, Tataru\cite{tataru} has shown that the above conjecture can be true in high 
dimension if the Sobolev space $\dot H^s$ is
replaced by the slightly smaller Besov space $\dot B^{s,1}_2$, which we define by
$$ \| f\|_{\dot B^{s,1}_2} = \sum_j 2^{js} (\int_{|\xi| \sim 2^j} |\hat{f}(\xi)|^2\ d\xi)^{1/2}.$$
More precisely, we have

\begin{theorem}\label{tat}\cite{tataru}  The Cauchy problem \eqref{wave-map} is well-posed for
small data in $\dot B^{\frac{n}{2},1}_2$ if $n > 3$.
\end{theorem}

In the case of the circle $m=2$, solutions to \eqref{wave-map} are complex exponentials of solutions
to the free wave equation, and both of the above results are 
true in any dimension $n$.  
However, when $m > 2$ the above Conjecture and Theorem do not hold in one dimension for $m>2$.  More precisely, we shall show

\begin{theorem}\label{main}  Let $n=1$ and $m>2$, and let $\eps > 0$, $C>0$.  Then there exists initial data $f$ with
\be{small}
\|f-e_1\|_{C^\infty_0} \lesssim \eps, \quad f-e_1 \hbox{ supported on } [-C,C]
\end{equation}
such that the solution  $\phi$ to \eqref{wave-map} satisfies
$$ \| \phi(T)\|_{\dot H^{1/2}} \to \infty \hbox{ as } T \to \infty.$$
In particular, we also have
$$ \| \phi(T)\|_{\dot B^{1/2,1}_2} \to \infty \hbox{ as } T \to \infty.$$
In fact, the norms grow like $(\log T)^{1/2}$ and $\log T$ respectively.  
\end{theorem}

By rescaling this theorem one may find
data with
$\| f \|_{\dot B^{1/2,1}_2} \lesssim \eps$ and
$\| \phi(\eps) \|_{\dot H^{1/2}}$ arbitrarily large.
Applying this fact with $\eps$ replaced
by $\eps 2^{-k}$ for $k=1,2,\ldots$ and piecing together different solutions via
finite speed of propagation we thus obtain

\begin{corollary}  Conjecture \ref{hs-c} and the claim of Theorem \ref{tat} are false
when $n=1$ and $m>2$.  In fact, one can construct data with arbitrarily small
norm in $\dot B^{1/2,1}_2$ (and thus in $\dot H^{1/2}$) which
cannot be continued for any non-zero time in $\dot B^{1/2,1}_2$ or even in 
$\dot H^{1/2}$.
\end{corollary}

We remark that even though the one-dimensional wave map problem
appears to be badly behaved in the spaces $\dot H^{1/2}$,
$\dot B^{1/2,1}_2$, one does have global existence and scattering
for large data in the critical space $L^{1,1}$ of functions with
derivatives in $L^1$.  See \cite{keel:wavemap}.

We now begin the proof of the theorem.  Let $f$ be a function obeying \eqref{small},
and let $\phi$ be the corresponding global solution to \eqref{wave-map}.
(For a proof that the wave map problem is globally well-posed for
smooth data in $n=1$, see \cite{gu}).  

To study the global behaviour of $\phi$, we introduce null co-ordinates
\begin{align*}
u & = x + t  \qquad\qquad v = x - t\\
\partial_u & = \frac{1}{2}(\partial_x + \partial_t)  \quad  \partial_v = \frac{1}{2}(\partial_x - \partial_t). \label{nullcoords}
\end{align*}
From \eqref{wave-map} we then have
\be{nul}
\phi_{uv} = - \phi (\phi_u \cdot \phi_v).
\end{equation}
Since $\phi$ stays on the sphere, we have $\phi_u \cdot \phi = \phi_v \cdot \phi = 0$.  From this and \eqref{nul} we see that
$\phi_u, \phi_v$ are both orthogonal to $\phi_{uv}$.  This implies the
following pointwise conservation laws, first observed by Pohlmeyer \cite{pohl}:
$$ \partial_v (|\phi_u|^2) = 2 \phi_u \cdot \phi_{uv} = 0, 
\quad \partial_u (|\phi_v|^2) = 2 \phi_v \cdot \phi_{uv} = 0.$$
One can also obtain these conservation laws from the conformal invariance of \eqref{wave-map}, or
from the trace-free property of the stress-energy tensor.  See e.g. \cite{shatah}, \cite{keel:wavemap}.
These laws are special to the one-dimensional case; in higher dimensions the wave map equation
is not completely integrable.

Thus $|\phi_u|$ is constant in the $v$ direction, and $|\phi_v|$ is constant in the $u$ direction.  Combining
this fact with the support conditions on the initial data, we thus obtain (after converting back into
spacetime co-ordinates)
\bas
(\partial_x + \partial_t) \phi(t,x) = 0 &\hbox{ when } |t+x|\geq C,\\
(\partial_x - \partial_t) \phi(t,x) = 0 &\hbox{ when } |t-x|\geq C.
\end{align*}
In particular, $\phi$ is constant on the regions $\{x < -|t|-C\}$, $\{x > |t|+C\}$, and $\{t > |x|+C\}$.
In the first two regions we have $\phi=e_1$ from the conditions on the initial data, while on the
last region we write $\phi=\alpha$ where $\alpha$ is a constant.  As we shall see, the question
of whether $\alpha = e_1$ will be crucial.

\begin{figure}[htbp] \centering
\ \psfig{figure=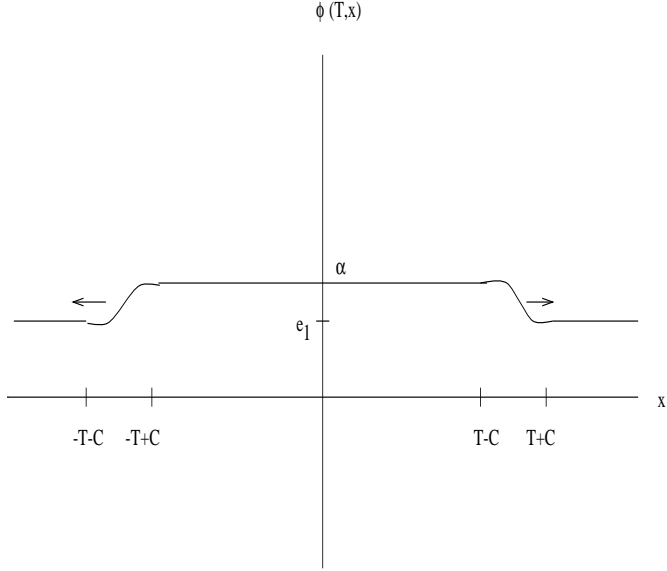,height=3in,width=3.6in}
\caption{A schematic depiction of $\phi(T,x)$ for some large $T$.  The
$\dot H^{1/2}$ norm grows logarithmically with $T$ if $\alpha \neq e_1$.
 }
        \label{fig}
        \end{figure}

In the regions $\{|t+x|< C, t > C\}$ and $\{|t-x|< C, t > C\}$, $\phi$ is a travelling wave
in the $\partial_x - \partial_t$ and $\partial_x + \partial_t$ directions respectively.  We
thus can describe $\phi(t,x)$ for $t > C$ by
\be{desc}
\phi(t,x) = 
\left\{\begin{array}{ll}
e_1& x \leq -t-C\\
F(t+x)& -t-C \leq x \leq -t+C\\
\alpha& -t+C \leq x \leq t-C\\
G(t-x)& t-C \leq x \leq t+C\\
e_1& t+C \leq x
\end{array}\right.
\end{equation}
where $F$, $G$ are $C^\infty$ functions on $[-C,C]$ such that $F(-C)=G(-C)=e_1$ and $F(C)=G(C)=\alpha$.  See Figure \ref{fig}.
We extend $F(x)$, $G(x)$ to equal $e_1$ for $x<-C$ and to equal $\alpha$ when $x>C$; note that 
$F$, $G$ remain smooth by \eqref{desc}
and the smoothness of $\phi$.

We can now compute the $\dot H^{1/2}$ norm of $\phi(T)$ for $T \gg C$.  Differentiating \eqref{desc} we have
$$ \phi_x(T,x) = F^\prime(T+x) - G^\prime(T-x)$$
and so
$$ \| \phi(T)\|_{\dot H^{1/2}} = \| F^\prime(T+\cdot) - G^\prime(T-\cdot) \|_{\dot H^{-1/2}}.$$
The spatial Fourier transform of $F^\prime(T+x) + G^\prime(T-x)$ is
$$ e^{-2\pi i T \xi} A(\xi) - e^{2\pi i T \xi} B(-\xi),$$
where $A$ and $B$ are the Fourier transforms of the Schwarz functions
$F^\prime$ and $G^\prime$ respectively.  We thus have
$$  \| \phi(T)\|_{\dot H^{1/2}}^2 =  \int \frac{|e^{-2\pi i T \xi} A(\xi) - e^{2\pi i T \xi} B(-\xi)|^2}{|\xi|}
\ d\xi.$$
Clearly this quantity is greater than
\be{test}
  \int_{T^{-1} \ll \xi \ll 1} 
\frac{|e^{-2\pi i T \xi} A(\xi) - e^{2\pi i T \xi} B(-\xi)|^2}{|\xi|}
\ d\xi.
\end{equation}
The functions $A$ and $B$ are Schwarz, and so we have the estimate $A(\xi) = A(0) + O(|\xi|)$,
$B(\xi) = B(0) + O(|\xi|)$.  In particular we have
$$|e^{-2\pi i T \xi} A(\xi) - e^{2\pi i T \xi} B(-\xi)|^2 = |e^{-2\pi i T \xi} A(0) - e^{2\pi i T \xi} B(0)|^2
+ O(|\xi|).$$
However, from the definition of $A$, $B$ and the boundary conditions of $F$, $G$ we see that $A(0)=B(0)=\alpha-e_1$,
so that
$$|e^{-2\pi i T \xi} A(\xi) - e^{2\pi i T \xi} B(-\xi)|^2 = 4 |\alpha - e_1|^2 \sin^2(2\pi T\xi) + O(|\xi|).$$
Inserting this into \eqref{test} we obtain
$$ \| \phi(T)\|_{\dot H^{1/2}}^2 \gtrsim |\alpha - e_1|^2  \int_{T^{-1} \ll \xi \ll 1} \frac{\sin^2(2 \pi T\xi)}{|\xi|}\ d\xi
+ O(1) \gtrsim |\alpha - e_1|^2 \log T + O(1).$$
This proves the Theorem, providing that we can construct $f$, $g$ such that
$$ \alpha \neq e_1.$$
In the case $m=1$ this is impossible since one has the explicit solution $\phi(t,x) = f(x+t)^{1/2} f(x-t)^{1/2}$
in complex co-ordinates.  However we shall show that one can make $\alpha \neq e_1$ when $m > 2$.  We shall do this by expanding $\alpha$ as a power series for data
close to $e_1$.  It turns out that one has to expand quite far in order to do
this, because the deviation of $\alpha$ is extremely
small.  In fact we have $\alpha = e_1 + O(\eps^5)$.

\section{The deviation of $\alpha$} 

Suppose for contradiction that $\alpha = e_1$ for all choices of data $f$ satisfying \eqref{small}.  In other words, we assume
$$ \phi(C,0) = e_1$$
for all such $\phi$.  Since our initial velocity is zero, the Cauchy problem
is time-symmetric and so $\phi(t,x) = \phi(-t,x)$ for all $t,x$.  Thus we also 
have $\phi(-C,0) = e_1$.  From \eqref{small} we have $\phi(0,\pm C) = e_1$.
By two applications of the fundamental theorem of calculus we therefore have
(in null co-ordinates)
$$ \int_{-C}^C \int_{-C}^C \phi_{uv}\ du dv = 0$$
for all $\phi$.  By \eqref{nul} this implies
\be{contra} \int_{-C}^C \int_{-C}^C \phi (\phi_u \cdot \phi_v)\ du dv = 0.
\end{equation}

Let $h$
be a bump function on $[-C,C]$ taking values in $\R^{m-1}$, the subspace of $\R^m$ 
orthogonal to $e_1$; the exact choice of $h$ will be determined later.  We define the initial data
$f(x)$ to be
\be{f-def} f(x) = \sum_{n=0}^\infty \frac{\eps^n}{n!} h(x)^n
= e_1 + \eps h(x) - \frac{\eps^2}{2} |h(x)|^2 e_1 - \frac{\eps^3}{6} |h(x)|^2 h(x) + \ldots
\end{equation}
where the expression $h(x)^n$ is defined by
$$ h(x)^n = 
\left\{\begin{array}{ll}
(-1)^{n/2} |h(x)|^n e_1 & \hbox{ when } n \hbox{ is even}\\
(-1)^{(n-1)/2} |h(x)|^{n-1} h(x) & \hbox{ when } n \hbox{ is odd.}
\end{array}\right.
$$
Note that $f$ takes values on the sphere and obeys \eqref{small}.

By the classical well-posedness theory we may obtain the asymptotic expansion
\be{asym} \phi = \phi^0 + \eps \phi^1 + \eps^2 \phi^2 + \eps^3 \phi^3 + \eps^4 \phi^4
+ O(\eps^5)
\end{equation}
for times $|T| < 2C$, if $\eps$ is sufficiently small.  Here $\phi^i$ are smooth
functions independent of $\eps$ for $i = 0, \ldots, 4$, and the error term
has a $C^{10}$ norm of $O(\eps^5)$.

Setting $\eps=0$, we see that $\phi^0 = e_1$ by the classical uniqueness theory.  
If one then inserts \eqref{asym}
into \eqref{nul} and extracts co-efficients, one obtains
\begin{align}
\phi^1_{uv} &= 0 \label{phi1-eq}\\
\phi^2_{uv} &= - e_1 (\phi^1_u \cdot \phi^1_v) \label{phi2-eq}\\
\phi^3_{uv} &= - e_1 (\phi^1_u \cdot \phi^2_v + \phi^2_u \cdot \phi^1_v)
- \phi^1 (\phi^1_u \cdot \phi^1_v) \label{phi3-eq}\\
\phi^4_{uv} &= - e_1 (\phi^1_u \cdot \phi^3_v + \phi^2_u \cdot \phi^2_v
+ \phi^3_u \cdot \phi^1_v) - \phi^1 (\phi^1_u \cdot \phi^2_v + \phi^2_u \cdot \phi^1_v)
- \phi^2 (\phi^1_u \cdot \phi^1_v)\label{phi4-eq}
\end{align}
Also, from \eqref{f-def} we have
\begin{align}
\phi^1(u,u) &= h(u) \label{phi1}\\
\phi^2(u,u) &= -\frac{1}{2}|h(u)|^2 e_1\label{phi2}\\
\phi^3(u,u) &= -\frac{1}{6}|h(u)|^2 h(u)\label{phi3}\\
\phi^4(u,u) &= \frac{1}{24} |h(u)|^4 e_1\label{phi4}
\end{align}
Finally, from time-reversal symmetry we have $\phi^i(u,v) = \phi^i(v,u)$ for
all $i = 1, \ldots, 4$ and all $u,v$.  We can now proceed to compute
$\phi^1, \phi^2, \ldots$ iteratively.

One can see inductively that $\phi_i$ is orthogonal to $e_1$ when $i$ is odd,
and parallel to $e_1$ when $i$ is even, for $i = 1, \ldots, 4$.  In particular,
the $(\phi^1_u \cdot \phi^2_v + \phi^2_u \cdot \phi^1_v)$ terms in
\eqref{phi3-eq} and \eqref{phi4-eq} vanish.

Solving the time-symmetric Cauchy problem \eqref{phi1-eq} with initial
data \eqref{phi1} we obtain
\be{phi1-form}
\phi^1(u,v) = \frac{h(u) + h(v)}{2}.
\end{equation}
Substituting \eqref{phi1-form}
into \eqref{phi2-eq} and solving the time-symmetric Cauchy problem with initial
data \eqref{phi2} we obtain
\be{phi2-form}
 \phi^2(u,v) = -\frac{1}{8} |h(u) + h(v)|^2 e_1.
\end{equation}

If we substitute \eqref{phi1-form} into \eqref{phi3-eq}
we obtain
\be{phi3-uv2}
\phi^3_{uv}(u,v) = - \frac{1}{4} \phi^1(u,v) (h'(u) \cdot h'(v)).
\end{equation}

We observe the following property of $\phi^1 \cdot \phi^3$:

\begin{lemma}\label{record}  We have $(\phi^1 \cdot \phi^3)(u,C) = 
(\phi^1 \cdot \phi^3)(u,-C)$ for all $u$.
\end{lemma}

\begin{proof}  From \eqref{phi1-form} and the fact that $h(\pm C) = 0$, we can rewrite this as
\be{h-targ}
 h(u) \cdot \phi^3(u,C) - h(u) \cdot \phi^3(u,-C) = 0.
\end{equation}

We now invoke the explicit formula
\be{phi3-old} \phi^3(u,v) = -\frac{|h(u)|^2 h(u) + |h(v)|^2 h(v)}{12}
+ (H(u) - H(v)) (h(u) - h(v))
\end{equation}
for $\phi^3$, where the matrix-valued function $H(u)$ is defined by
$$ H'(u) = \frac{1}{8} h(u) h'(u)^t, \quad H(-C) = 0.$$
Indeed, one easily verifies that the proposed solution $\phi^3$ in \eqref{phi3-old} satisfies \eqref{phi3}, is time symmetric, and satisfies
$$ \phi^3_{uv} = -H'(u)h'(v) - H'(v)h'(u) = - \frac{1}{8} h(u) (h'(u) \cdot h'(v))
- \frac{1}{8} h(v) (h'(v) \cdot h'(u))$$
which is \eqref{phi3-uv2}. 

Substituting \eqref{phi3-old} into \eqref{h-targ} and noting that $h(\pm C), H(-C) = 0$, we reduce to
$$ h(u) \cdot H(C) h(u) = 0.$$
It thus suffices to show that $H(C)$ is an anti-symmetric matrix.  To see this,
we use  the fundamental theorem of calculus to obtain
$$ H(C) = \int_{-C}^C \frac{1}{8} h(u) h'(u)^t\ du.$$
The symmetric part of this is
$$ \frac{H(C) + H(C)^t}{2} = 
\int_{-C}^C \frac{1}{16} (h(u) h'(u)^t + h'(u) h(u)^t)\ du
= \frac{1}{16} h(u) h(u)^t \bigr|_{-C}^C = 0$$
as desired.
\end{proof}

We now expand \eqref{contra} as a power series in $\eps$.  Remarkably,
one must go up to the $\eps^5$ term in order to get an expression
which does not automatically vanish.  
Since the derivatives of $\phi_0, \phi_2, \phi_4$ are orthogonal to the
derivatives of $\phi_1$, $\phi_3$,
the $\eps^5$ component of \eqref{contra} is
\be{contra2} \int_{-C}^C \int_{-C}^C 
\phi^1 (\phi^1_u \cdot \phi^3_v + \phi^2_u \cdot \phi^2_v + \phi^3_u \cdot \phi^1_v)
+ \phi^3 (\phi^1_u \cdot \phi^1_v)\ du dv = 0.
\end{equation}
The idea is to choose $h$ so that the left-hand side of \eqref{contra2}
is non-zero, thus obtaining the desired contradiction.  Although one could
in principle achieve this by applying \eqref{phi3-old}, we shall
instead attempt to simplify \eqref{contra2} by
repeated integrations by parts.

From the identity
$$ \phi^1_u \cdot \phi^3_v + \phi^3_u \cdot \phi^1_v
= (\phi^1 \cdot \phi^3)_{uv} - \phi^1 \cdot \phi^3_{uv}$$
which follow from the product rule and \eqref{phi1-eq}, we can
rewrite \eqref{contra2} as
\be{multi} \int_{-C}^C \int_{-C}^C \phi^1 (\phi^1 \cdot \phi^3)_{uv} +
\phi^1(\phi^2_u \cdot \phi^2_v - \phi^1 \cdot \phi^3_{uv})
+ \phi^3 (\phi^1_u \cdot \phi^1_v)\ du dv = 0.
\end{equation}
We now claim that
$$ \int_{-C}^C \int_{-C}^C \phi^1 (\phi^1 \cdot \phi^3)_{uv}\ du dv = 0.$$
To see this, it suffices by \eqref{phi1-form} and symmetry to show that
$$ \int_{-C}^C \int_{-C}^C h(u) (\phi^1 \cdot \phi^3)_{uv}\ du dv = 0.$$
Integrating in the $v$ variable we reduce to
$$ \int_{-C}^C h(u) (\phi^1 \cdot \phi^3)_{u}|_{v=-C}^{v=C}\ du = 0.$$
But this follows from differentiating Lemma \ref{record} with respect to $u$.

Next, we treat the integral
$$ \int_{-C}^C \int_{-C}^C \phi^3 (\phi^1_u \cdot \phi^1_v)\ du dv.$$
By \eqref{phi1-form}, this is
$$ \frac{1}{4} \int_{-C}^C \int_{-C}^C \phi^3 (h'(u) \cdot h'(v))\ du dv.$$
Integrating by parts twice and using the fact that $h(\pm C) = 0$, this is
$$ \frac{1}{4} \int_{-C}^C \int_{-C}^C \phi^3_{uv} (h(u) \cdot h(v))\ du dv.$$
Inserting these identities into \eqref{multi}, we reduce to
$$
\int_{-C}^C \int_{-C}^C 
\phi^1(\phi^2_u \cdot \phi^2_v - \phi^1 \cdot \phi^3_{uv})
+ \frac{1}{4} \phi^3_{uv} (h(u) \cdot h(v))\ du dv = 0.
$$
By \eqref{phi3-uv2}, this becomes
$$
\int_{-C}^C \int_{-C}^C 
\phi^1[\phi^2_u \cdot \phi^2_v + \frac{|\phi^1|^2}{4} h'(u) \cdot h'(v)	
- \frac{h'(u) \cdot h'(v)}{16} h(u) \cdot h(v)]\ du dv = 0.$$
By \eqref{phi2-form} and \eqref{phi1-form}, this is
\bas
\int_{-C}^C \int_{-C}^C 
\phi^1[&\frac{(h(u) \cdot h'(u) + h(v) \cdot h'(u)}{4}
\ \frac{h(u) \cdot h'(v) + h(v) \cdot h'(v)}{4}\\
&+ \frac{|h(u)+h(v)|^2}{16} h'(u) \cdot h'(v)
-\frac{h'(u) \cdot h'(v)}{16} h(u) \cdot h(v)]\ du dv = 0.
\end{align*}
By symmetry we may replace the $\phi^1$ with an $h(u)$.
Multiplying by 16 and expanding, we thus obtain
\bas
\int_{-C}^C \int_{-C}^C
h(u)[&(h(u) \cdot h'(u)) (h(u) \cdot h'(v)) + (h(u) \cdot h'(u)) (h(v) \cdot h'(v))\\
& + (h(v) \cdot h'(u)) (h(u) \cdot h'(v)) +
(h(v) \cdot h'(u)) (h(v) \cdot h'(v))\\
&+ (|h(u)|^2 + h(u) \cdot h(v) + |h(v)|^2) h'(u) \cdot h'(v)]\ du dv = 0.
\end{align*}
The first and second terms vanish 
since $h'(v)$ and $h(v) \cdot h'(v)$ both have mean zero.
We thus reduce to
\bas
\int_{-C}^C \int_{-C}^C
h(u)[&(h(v) \cdot h'(u)) (h(u) \cdot h'(v)) +
(h(v) \cdot h'(u)) (h(v) \cdot h'(v))\\
&+ (|h(u)|^2 + h(u) \cdot h(v) + |h(v)|^2) h'(u) \cdot h'(v)]\ du dv = 0.
\end{align*}
We now utilize the hypothesis $m>2$, and
choose $h(x) = h_2(x) e_2 + h_3(x) e_3$.
Extracting the $e_2$ component of the above equation, we obtain
\be{big}
\bs
\int_{-C}^C \int_{-C}^C &h_2(u)[ (h_2(v) h_2'(u) + h_3(v) h_3'(u))(h_2(u) h_2'(v) + h_3(u) h_3'(v))\\
&+ (h_2(v) h_2'(u) + h_3(v) h_3'(u)) (h_2(v) h_2'(v) + h_3(v) h_3'(v))\\
&+ (h_2(u)^2 + h_3(u)^2 + h_2(u) h_2(v) + h_3(u) h_3(v) + h_2(v)^2 + h_3(v)^2)\\
&(h_2'(u) h_2'(v) + h_3'(u) h_3'(v))]\ du dv = 0.
\end{split}
\end{equation}
The integrand splits into several expressions which are the
product of a function of $u$ and a function of $v$.  If one of these functions has mean zero then its contribution to \eqref{big} vanishes.  Examples of
such functions include $h_2'$, $h_3'$, $h_2 h_2'$, $h_3 h_3'$, $h_2^2 h_2'$, $h_3^2 h_3'$.  Eliminating all such terms from \eqref{big}, one is reduced to
\bas
\int_{-C}^C \int_{-C}^C
h_2(u)[ &h_2(v) h_2'(u) h_3(u) h_3'(v)
 + h_3(v) h_3'(u) h_2(u) h_2'(v) \\
&+ h_3(v) h_3'(u) h_2(v) h_2'(v)
+ h_2(u) h_2(v) h_3'(u) h_3'(v)\\
&+ h_3(u) h_3(v) h_2'(u) h_2'(v) + h_2(v)^2 h_3'(u) h_3'(v) ]\ du dv = 0.
\end{align*}
We split the $u$ and $v$ integrations to obtain
\be{da} DA + EB + AD + EA + DB + AE = 0
\end{equation}
where 
\bas
A &= \int_{-C}^C h_2(v) h_3'(v)\ dv\\
B &= \int_{-C}^C h_3(v) h_2'(v)\ dv\\
D &= \int_{-C}^C h_2(u) h_2'(u) h_3(u)\ du\\
E &= \int_{-C}^C h_2(u)^2 h_3'(u)\ du
\end{align*}
By integrating by parts we see that $B=-A$ and $D=-E/2$.  The left-hand side
of \eqref{da} therefore simplifies to $AE/2$.  To obtain a contradiction
it thus suffices to choose $h_2$, $h_3$ so that $A,E \neq 0$.
But this is easily achieved, e.g. let $h_3$ be such that
$\int_{-C}^C h_3'(u)^3\ du \neq 0$, and set $h_2 = h_3'$.

This shows that for any $\eps > 0$ one can find initial data
for which $\alpha \neq e_1$, which concludes the proof
of Theorem \ref{main}.

\endprf

\section{Further Remarks}

The ill-posedness result is not confined to the sphere $S^2$, and can
in fact be generalized to other manifolds with positive or negative
curvature.  For the general version of the Pohlmeyer identity used
in the above argument, see \cite{keel:wavemap}.  The fact that
$\alpha \neq e_1$ for general manifolds with non-zero Riemmanian
curvature follows by a tedious but straightforward expansion of
$\phi$ as was done above.

The wave map equation \eqref{wavemap} can be written in integral form as
$$ \phi = \phi_0 + \frac{1}{2} \Box^{-1} (\phi \Box (\phi \cdot \phi))
- \frac{1}{2} \Box^{-1} (\phi (\phi \cdot \Box \phi))$$
where $\phi_0$ denotes the solution to the free wave equation with
initial position $f$ and zero initial velocity, and $\Box^{-1}$
is the fundamental solution to the inhomogeneous Cauchy problem for
the linear wave equation.  

The standard way to show well-posedness of these Cauchy problems
is to iterate the above integral equation and hope that the
solution converges in the desired norm.  In the $1+1$
dimensional problem, however, if one starts with the free solution
$\phi_0 = \frac{1}{2}(f(u) + f(v))$ and computes a few 
iterates, one is soon led to expressions such as $D^{-1}(g Dh) (u)$,
where $g(u)$ and $h(u)$ are one-dimensional functions which have
the same regularity as the initial data $f$, and $D = \frac{d}{du}$
is the usual differentiation operator.  These expressions were already
seen in the previous section.

For subcritical regularities, these types of expressions are
well-behaved if we localize in the $u$ variable.  However at the
critical regularity it is pointless to localize (as can be seen by scale
invariance considerations), and
the bad behaviour of the expression $D^{-1}(g Dh)$ then dooms
any attempt to obtain well-posedness by iterative methods.  
Even if $h$ was a $C^{\infty}_0$
function, we could take $g=Dh$, so that $g Dh$ would be non-negative,
and $D^{-1}(g Dh)$ would be a smoothed out Heaviside function,
which barely fails to be in the space $B^{1/2,2}_1$ or $\dot H^{1/2}$.
Indeed, the Fourier transform of such a function behaves like $1/\xi$ near
the origin.  On the other hand, this function is in the space $L^{1,1}$
of functions with absolutely integrable derivative, which is compatible
with the results in \cite{keel:wavemap}.

The expression $\Box^{-1}(\phi \Box \psi)$ appears to be better behaved 
when the dimension $n$ is large, so this difficulty may well be isolated
to the one-dimensional case.  One might also argue that this ill-posedness
result is somehow related to the pathologies of the negative order Sobolev
space $\dot H^{-1/2}$, and that such spaces do not occur in the 
higher-dimensional theory.

Daniel Tataru has informed us that Theorem \ref{tat} can be
extended to $n > 1$.  Thus there is a sharp distinction between
the one-dimensional and higher-dimensional cases.

\section{Acknowledgements}

The author thanks Sergiu Klainerman for his hospitality, encouragement, and
very informative discussions while this work was performed.  The author
is also deeply indebted to the reviewer for pointing out an error in the
first version of this paper.
The author is partially supported
by NSF grant DMS-9706764.


\begin{thebibliography}{10}

\bibitem{grillakis.zurich}
M. Grillakis, \emph{A priori estimates and regularity of nonlinear waves},
in Proceed. Inter. Congress of Math. 1994, Birkh\"auser, 1187 - 1194.

\bibitem{gu}
C. Gu, \emph{On the Cauchy problem for harmonic maps defined on two-dimensional
Minkowski space},
Comm. Pure Appl. Math., \textbf{33},(1980),
727--737.

\bibitem{keel:wavemap}
M. Keel, T. Tao, \emph{Local and global well-posedness
of wave maps on $\R^{1+1}$ for rough data}, IMRN \textbf{21} (1998),
1117--1156.

\bibitem{kman.barrett}
S. Klainerman, \emph{On the regularity of classical field theories in Minkowski
space-time $\R^{3+1}$}, Prog. in Nonlin. Diff. Eq. and their Applic., \textbf{
29},
(1997), Birkh\"auser, 113--150.

\bibitem{kman.mach.smoothing}
S. Klainerman, M. Machedon, \emph{Smoothing estimates for null forms
and applications}, Duke Math. J., \textbf{81} (1995), 99--133.
 
\bibitem{kman.selberg}
S. Klainerman, S. Selberg, \emph{ Remark on the optimal regularity for
equations of wave maps type}, C.P.D.E., \textbf{22} (1997), 901--918.

\bibitem{lady}
O.A. Ladyzhenskaya, V.I. Shubov, \emph{Unique solvability of the Cauchy
problem for the equations of the two dimensional chiral fields, taking values
in complete Riemann manifolds}, J. Soviet Math.,
\textbf{25} (1984), 855--864. (English Trans. of 1981 Article.)
 
\bibitem{pohl}
K. Pohlmeyer, \emph{Integrable Hamiltonian systems and interaction through
quadratic constraints}, Comm. Math. Phys.,
\textbf{46} (1976),
207--221.

\bibitem{shatah}
J. Shatah, \emph{Weak solutions and development of singularities
of the $SU(2)$ $\sigma$-model.}
Comm. Pure Appl. Math., \textbf{41} (1988),
459--469.

\bibitem{shatah.zurich}
J. Shatah, \emph{The Cauchy problem for harmonic maps on Minkowski space},
in Proceed. Inter. Congress of Math. 1994, Birkh\"auser, 1126--1132.

\bibitem{shatah.shadi.blow} 
ow
J. Shatah, A. Tavildar-Zadeh, \emph{On the Cauchy problem for equivariant
wave maps}, Comm. Pure Appl. Math., \textbf{47} (1994), 719 - 753.

\bibitem{struwe.barrett}
M. Struwe, \emph{Wave Maps, in Nonlinear Partial Differential Equations
in Geometry and Physics}, Prog. in Nonlin. Diff. Eq. and their Applic., 
\textbf{29}, (1997), Birkh\"auser, 113--150.

\bibitem{tataru}
D. Tataru, \emph{Local and global results for wave maps I},
Preprint, 1997.

\end{thebibliography}
\end{document}